\documentclass[12pt,a4paper]{article}
\usepackage{amsmath}
\usepackage{amssymb}
\usepackage{hanging}

\begin{document}
\begin{center}
\Large{On reducing inconsistency of pairwise comparison matrices below an acceptance threshold}\\[1cm]
\end{center}
\begin{center}
\normalsize
S\'andor BOZ\'OKI$^{1}$ \\
\small{
Laboratory on Engineering and Management Intelligence, \\
Research Group of Operations Research and Decision Systems, \\
Institute for Computer Science and Control, \\
Hungarian Academy of Sciences (MTA SZTAKI);
and \\
Department of Operations Research and Actuarial Sciences \\
Corvinus University of Budapest, Hungary \\
\verb|bozoki.sandor@sztaki.mta.hu|  } \\[5mm]
\end{center}
\begin{center}
\normalsize
\normalsize
J\'anos F\"UL\"OP \\
\small{
Laboratory on Engineering and Management Intelligence, \\
Research Group of Operations Research and Decision Systems, \\
Institute for Computer Science and Control, \\
Hungarian Academy of Sciences (MTA SZTAKI);
Budapest, Hungary \\
\verb|fulop.janos@sztaki.mta.hu|  } \\[5mm]
\end{center}
\begin{center}
\normalsize
Attila POESZ \\
\small{
Department of Operations Research and Actuarial Sciences \\
Corvinus University of Budapest, Hungary \\
\verb|attila.poesz@uni-corvinus.hu|  } \\[25mm]
\end{center}
\begin{center}
\normalsize{4 November 2013}\\[4cm]
\end{center}
\begin{center}
Research was supported in part by OTKA grant K77420.
\end{center}
\footnotetext[1]{corresponding author}
\newpage

\begin{abstract}
A recent work of the authors on the analysis of
pairwise comparison matrices that can be made consistent by the modif{\kern0pt}ication
of a few elements is continued and extended.
Inconsistency indices  are def{\kern0pt}ined for indicating the
overall quality of a pairwise comparison matrix.
It is expected that serious contradictions in the matrix
imply high inconsistency and vice versa.
However, in the 35-year history of the
applications of pairwise comparison matrices, only one of
the indices, namely $CR$ proposed by Saaty, has been associated
to a general level of acceptance, by the well known ten percent rule.
In the paper, we consider a wide class of inconsistency indices,
including $CR$, $CM$ proposed by Koczkodaj and Duszak
and $CI$ by Pel\'aez and Lamata.
Assume that a threshold of acceptable inconsistency is given
(for $CR$ it can be 0.1).
The aim is to f{\kern0pt}ind the minimal number of matrix elements,
the appropriate modif{\kern0pt}ication of which makes the matrix
acceptable.
On the other hand, given the maximal number of
modif{\kern0pt}iable matrix elements, the aim is to f{\kern0pt}ind
the minimal level of inconsistency that can be achieved.
In both cases the solution is derived from a nonlinear mixed-integer
optimization problem. Results are applicable in
decision support systems that allow real time interaction
with the decision maker in order to review pairwise
comparison matrices.
\end{abstract}

\section{Introduction}
\label{intro}

Pairwise comparison matrices (Saaty, 1977) are used in multi-attribute
decision problems, where relative importance of the criteria,
the evaluations of the alternatives with respect to each
criterion are to be quantif{\kern0pt}ied. The method of pairwise comparison
is also applied for determining voting powers
in group decision making. One of the advantages of pairwise
comparison matrices is that the decision maker is faced to
a sequence of elementary questions concerning the comparison of two
criteria/alternatives at a time, instead of a complex task
of providing the weights of the whole set of them.

A real $n \times n$ matrix  $A$ is a \emph{pairwise comparison matrix} if it is positive and reciprocal, i.e.,\
\begin{eqnarray}
a_{ij}&>&0, \label{eq:positive}\\
a_{ij}&=&\frac{1}{a_{ji}} \label{eq:reciprocal}
\end{eqnarray}
for all $i,j=1,\dots,n$. ${A}$ is \emph{consistent} if the transitivity property
\begin{align}
a_{ij} a_{jk} = a_{ik}   \label{eq:TransitivityProperty}   %(\ref{eq:TransitivityProperty})
\end{align}
holds for all $i,j,k=1,2,\ldots,n$; otherwise it is called \emph{inconsistent}.

For a positive $n\times n$ matrix $A$, let $\bar A=\log A$ denote the $n\times n$ matrix with the elements
$$\bar a_{ij}=\log a_{ij}, \ \ i,j=1,\dots ,n.$$ Then $A$ is consistent if and only if
\begin{align}
\bar a_{ij}+\bar a_{jk}+\bar a_{ki}=0, \  \forall \, i,j,k=1,\dots ,n \label{eq:subspace}
\end{align}
holds. Matrices $\bar A$ fulf{\kern0pt}illing the homogenous linear system (\ref{eq:subspace})
constitute a linear subspace in  $\mathbb{R}^{n\times n}$.

Let ${\cal P}_n$ denote the set of the $n\times n$ pairwise comparison matrices, and
${\cal C}_n\subset{\cal P}_n$ the set of the consistent matrices.
Since the reciprocity constraint (\ref{eq:reciprocal}) corresponds to
$\bar a_{ij}=-\bar a_{ji}$
in the logarithmized space, the set $\log {\cal P}_n=\{\log A \mid A \in {\cal P}_n\}$ is the set of $n\times n$ skew-symmetric matrices,
an $n(n-1)/2$-dimensional linear subspace of $\mathbb{R}^{n\times n}$. The set $\log
{\cal C}_n=\{\log A \mid A\in {\cal C}_n\}$ is the set of matrices fulf{\kern0pt}illing (\ref{eq:subspace}),
and as pointed out in Chu (1997), is an $(n-1)$-dimensional linear subspace of $\mathbb{R}^{n\times n}$.
Clearly, $\log {\cal C}_n\subset \log {\cal P}_n$.

In decision problems of real life, the pairwise comparison matrices are rarely consistent.
Nevertheless, decision makers are interested in the level
of inconsistency of their judgements, which somehow expresses the goodness or
``quality'' of pairwise comparisons totally, because conf{\kern0pt}licting judgements may
lead to senseless decisions. Therefore, some index is needed to measure the possible
contradictions and inconsistencies of the pairwise comparison matrix.

A function $\phi_n : {\cal P}_n \to R$ is called an \emph{inconsistency index} if  $\phi_n (A)=0$ for every consistent and
$\phi_n (A)>0$ for every inconsistent pairwise comparison matrix $A$. The inconsistency indices used in the practice are continuous, and the value of $\phi_n (A)>0$ indicates, more or less, how much an inconsistent matrix dif{\kern0pt}fers from a consistent one.

Since in the practice the consistency of a pairwise comparison matrix is not easy to assure, certain level of inconsistency is usually accepted by the decision makers. This works in the practice in such a way that
for a given inconsistency index $\phi_n$ an acceptance threshold $\alpha_n \ge 0$ is chosen, and a
matrix $A\in {\cal P}_n$ is kept for further use only if $\phi_n (A)\le \alpha_n$ holds; otherwise, it is rejected or the pairwise comparisons are carried out again. The carrying out of all pairwise comparisons
for f{\kern0pt}illing-in the matrix is often a time-consuming task. Therefore, before the total rejection of a pairwise comparison matrix with an inconsistency level above a prescribes acceptance threshold, it may be worth investigating whether it is possible to improve the inconsistency of the matrix to an acceptable level by performing fewer pairwise comparisons.

The paper will concentrate on the following problem: for a given \linebreak  $A\in {\cal P}_n$, inconsistency index $\phi_n$ and acceptance level $\alpha_n$, what is the minimal number of the elements of matrix $A$ that by modifying these elements, and of course their reciprocals, the pairwise comparison matrix can be made acceptable. We shall show that under a slight boundedness assumption,
this can be achieved by solving a nonlinear mixed 0-1 optimization problem.  If it comes out that
the matrix can be turned into an acceptable one by modifying relatively few elements, then it may be a case when
a more-or-less consistent evaluator was less attentive at these few elements, or a data-recording error happened.
So it may be worth re-evaluating these elements. Of course, if the the evaluator insists on the previous  values,
or the acceptable inconsistency threshold cannot be reached with the new values, then this approach was unsuccessful: all pairwise comparisons are to be evaluated again. If however after the revision of the critical elements, the inconsistency level of the modif{\kern0pt}ied matrix is already acceptable, then we can continue the decision process with it.

Concerning the investigations above, when solving the nonlinear mixed 0-1 programming problems, it is very benef{\kern0pt}icial if the nonlinear optimization problems obtained after the relaxation of the 0-1 variables are convex optimization problems. In the convex case several sophisticated methods and softwares are available, while in the nonconvex case methodological and implementation dif{\kern0pt}f{\kern0pt}iculties may arise.  Since $\log {\cal C}_n$ is a linear subspace, ${\cal C}_n$ is a nonconvex manifold in $\mathbb{R}^{n\times n}$. One can immediately  conclude that it is better to investigate the convexity issues in the logarithmized space.

Several proposals of inconsistency indices are known, see the overviews of Brunelli and Fedrizzi (2011, 2013a)
and Brunelli et al. (2013b) for detailed lists and properties.
This paper focuses on three well-known inconsistency indices.
They are $CR$ proposed by Saaty (1980),
$CM$ proposed by Koczkodaj and Duszak (Koczkodaj 1993; Duszak and Koczkodaj 1994),
and $CI$ proposed by Pel\'aez and Lamata (2003).
The properties  and relationship of the fundamental indices  $CR$ and $CM$ were also studied in Boz\'oki and Rapcs\'ak (2008).
In this paper we point out that for the inconsistency indices in our focus, the nonlinear mixed 0-1 optimization problems mentioned above can be formulated in the logarithmized space, and appropriate convexity properties hold on them.
We show that $CR$ and $CI$ are convex function in the logarithmized space, and $CM$ is quasiconvex, but can be transformed into a convex function by applying a suitable strictly monotone univariate function on it.

This paper is in a close relation to an earlier paper of the authors (Boz\'oki et al.\ 2011b). In the latter paper we investigated the special case when the acceptance threshold $\alpha_n$ is 0, i.e.\ the modif{\kern0pt}ied pairwise comparison matrix must be consistent. No inconsistency indices were needed for this investigation, simple graph theoretic ideas were applied. Unfortunately, the technique applied for  $\alpha_n=0$ cannot be extended to the general case, therefore, a new approach is proposed in this paper.

We also mention that some of the issues investigated in this paper were already considered, in Hungarian, in Boz\'oki et al.\ (2012).

Since inconsistent matrices are in the focus of this paper, and for $n=1$ and $n=2$ the pairwise comparison matrices are consistent, we shall assume in the sequel, without loss of generality, that $n\ge 3$.

In Section 2, the optimization problems to be solved are presented in a general form. The general issues are specialized and investigated for the inconsistency indices $CR$ of Saaty, $CM$ of Koczkodaj and Duszak, and $CI$ proposed by Pel\'aez and Lamata in Sections 3 through 5, respectively. A numerical example is presented in Section 6.

\section{The general form of the optimization problems}

Let $\phi_n$ be an inconsistency index and $\alpha_n$ be an acceptance threshold, and let
\begin{align}
{\cal A}_n(\phi_n,\alpha_n)=\{A\in {\cal P}_n  \mid \phi_n (A)\le \alpha_n\}  \label{eq:Cal_A_n}
\end{align}
denote the set of $n\times n$ pairwise comparison matrices with
inconsistency $\phi_n$ not exceeding threshold $\alpha_n$.
Let $A, \hat A \in {\cal P}_n$ and
\begin{align}
d(A,\hat A)=\,\mid\!\{(i,j): 1\le i<j\le n, a_{ij}\ne \hat
a_{ij}\}\!\mid\label{eq:d}
\end{align}
denote the number of matrix elements above the main diagonal,
where matrices $A$ and $\hat A$ dif{\kern0pt}fer from each other.
By reciprocity, the number of dif{\kern0pt}ferent elements is
the same as in positions below the main diagonal.\\

Consider  pairwise comparison matrix  $A \in {\cal P}_n$ with
$\phi_n (A)> \alpha_n$ as it is not acceptable in terms of inconsistency.
We want to calculate the minimal number of matrix elements
above the main diagonal to be modif{\kern0pt}ied in order to make matrix acceptable
(elements below the main diagonal are determined by the
elements above the main diagonal). That is to solve
the optimization problem
\begin{equation}
\begin{array}{ll}
\min &d(A,\hat A)\\
{\rm s.t.} & \hat A\in {\cal A}_n(\phi_n,\alpha_n),
\end{array}
\label{eq:min_d}
\end{equation}
where the elements above the main diagonal of $\hat A$ are variables.

We could also ask the minimal inconsistency of  $A \in {\cal P}_n$
matrix can be reached by modifying at most $K$ elements and their reciprocals.
The optimization problem is
\begin{equation}
\begin{array}{ll}
\min &\alpha\\
{\rm s.t.} & d(A,\hat A)\le K,\\
 &\hat A\in {\cal A}_n(\phi_n,\alpha),
\end{array}
\label{eq:min_alpha}
\end{equation}
where $\alpha$ and the
elements above the main diagonal of $\hat A$ are variables.

Problems (\ref{eq:min_d}) and (\ref{eq:min_alpha}) can be formulated in logarithmic space:
\begin{align}
\log {\cal A}_n(\phi_n,\alpha_n)=\{X\in \log {\cal P}_n  \mid \phi_n (\exp
X)\le \alpha_n\}, \label{eq:log_Cal_A_n}
\end{align}
therefore (\ref{eq:min_d}) is equivalent to
\begin{equation}
\begin{array}{ll}
\min &d(\log A, X)\\
{\rm s.t.} & X\in \log {\cal P}_n,\\
& \phi_n(\exp X)\le \alpha_n,
\end{array}
\label{eq:log_min_d}
\end{equation}
where elements above the main diagonal of $X$ are variables.
The f{\kern0pt}irst constraint in (\ref{eq:log_min_d}) means that $X$ belongs to
the subspace of skew-symmetric matrices.
In this paper we show that the second, nonlinear inequality is a convex constraint
in case of inconsistency indices $CR$ (Saaty 1980), $CM$  (Koczkodaj 1993; Duszak and Koczkodaj 1994)
and $CI$ (Pel\'aez and Lamata, 2003).

Problem (\ref{eq:min_alpha}) can be rewritten in the same way as above:
\begin{equation}
\begin{array}{ll}
\min &\alpha\\
{\rm s.t.} & d(\log A,X)\le K,\\
  & X\in \log {\cal P}_n,\\
& \phi_n(\exp X )\le \alpha,
\end{array}
\label{eq:log_min_alpha}
\end{equation}
where $\alpha$  and elements above the main diagonal of $X$ are variables.

The objective function $d$ can be replaced by using the well-known \linebreak ``Big M'' method.
Assume  that $M\ge 1$ is given as an upper bound of the values of the
elements in  $A\in {\cal P}_n$ and the computed $\hat A \in {\cal P}_n$ matrices,
which is determined as the optimal solution of  problems (\ref{eq:min_d}) and (\ref{eq:min_alpha}), i.e.,
\begin{align}
1/M\le a_{ij}\le M,\ 1/M\le \hat a_{ij}\le M,\ i,j=1,\dots ,n.
\label{eq:bigm1}
\end{align}

We can f{\kern0pt}ind such an upper bound $M$ if we get a bounded interval by knowing
the actual level of $\phi_n$, which contains at least one optimal solution of
problems (\ref{eq:min_d}), and (\ref{eq:min_alpha}).

On the other hand, if a theoretical upper bound $M$ is not given,
then a reasonable bound $M$ is usually determined on the values of the pairwise
comparison matrices in every specif{\kern0pt}ic problem.
Constraint (\ref{eq:bigm1}) can be described as
\begin{align}
A, \hat A\in [1/M,M]^{n\times n} \label{eq:bigm2}
\end{align}
in matrix form, and if the condition (\ref{eq:bigm2}) associated
with $\hat A$ is attached to problems (\ref{eq:min_d}) and also (\ref{eq:min_alpha}), we get
\begin{equation}
\begin{array}{ll}
\min &d(A,\hat A)\\
{\rm s.t.} & \hat A\in {\cal A}_n(\phi_n,\alpha_n)\cap [1/M,M]^{n\times n},
\end{array}
\label{eq:min_d2}
\end{equation}
and, respectively,
\begin{equation}
\begin{array}{ll}
\min &\alpha\\
{\rm s.t.} & d(A,\hat A)\le K,\\
 &\hat A\in {\cal A}_n(\phi_n,\alpha)\cap [1/M,M]^{n\times n}.
\end{array}
\label{eq:min_alpha2}
\end{equation}

Introduce $\bar M=\log M$, problems (\ref{eq:min_d2}) and (\ref{eq:min_alpha2}) become equivalent to
\begin{equation}
\begin{array}{ll}
\min &d(\log A, X)\\
{\rm s.t.} & X\in \log {\cal P}_n \cap [-\bar M,\bar M]^{n\times n},\\
& \phi_n(\exp X)\le \alpha_n,
\end{array}
\label{eq:log_min_d2}
\end{equation}
and
\begin{equation}
\begin{array}{ll}
\min &\alpha\\
{\rm s.t.} & d(\log A,X)\le K,\\
  & X\in \log {\cal P}_n \cap [-\bar M,\bar M]^{n\times n},\\
& \phi_n(\exp X )\le \alpha.
\end{array}
\label{eq:log_min_alpha2}
\end{equation}
in the logarithmic space.

The ``Big M'' method can be applied for (\ref{eq:log_min_d2}) and (\ref{eq:log_min_alpha2}).
Let $\bar A=\log A$, and introduce binary variables $y_{ij}\in \{0,1\},\ 1\le i <j \le n$.
Problem (\ref{eq:log_min_d2}) can be altered by using $\bar A \in [-\bar M,\bar M]^{n\times n}$ into the following
mixed 0-1 programming problem:
    \begin{equation}
        \begin{array}{rllllll}
        \min         &&\sum\limits_{i=1}^{n-1} \sum\limits_{j=i+1}^{n} y_{ij}\\
        \mathrm{s.t.}&& \phi_n(\exp X)\le \alpha_n, \\
        && x_{ij}=-x_{ji},&& 1\le i \le j \le n,\\
        && -\bar M\le x_{ij} \le  \bar M, &&1\le i <j \le n,\\
        && -2\bar My_{ij}\le x_{ij}-\bar a_{ij} \le  2\bar My_{ij}, &&1\le i <j \le n,\\
        && y_{ij} \in \{0,1\},&&1\le i <j \le n.\\
        \end{array}
    \label{eq:log_min_d3}
    \end{equation}

The optimal value of (\ref{eq:log_min_d3}) gives the minimal number of the matrix elements
above the main diagonal to be modif{\kern0pt}ied in order to achieve $\phi_n \leq \alpha_n.$
In the optimal solution, $y_{ij}=1$ indicates the matrix elements that
(and their reciprocal pairs) are modif{\kern0pt}ied, and $\exp x_{ij}$ gives a feasible value of these elements.

Problem (\ref{eq:log_min_d3}) may have multiple optimal solutions
with respect to the binary variables. If all of them are of interest,
we list them one by one as follows.
Assume that $L^*$ is the optimum value of the problem (\ref{eq:log_min_d3}), \linebreak
$y_{ij}^*$,  $1\le i<j\le n$, is an optimal solution and
$I_0^*=\{(i,j)\mid y_{ij}^*=0, 1\le i<j\le n\}$. By adding the constraint
\begin{align}
\sum\limits_{i=1}^{n-1} \sum\limits_{j=i+1}^{n} y_{ij}=L^*
\label{eq:csat1}
\end{align}
 to (\ref{eq:log_min_d3}) we can ensure, that the optimal solutions of (\ref{eq:log_min_d3})
can only be the feasible solutions of (\ref{eq:log_min_d3})-(\ref{eq:csat1}).

The addition of constraint
\begin{align}
\sum\limits_{(i,j)\in I_0^*} y_{ij}\ge 1
\label{eq:csat2}
\end{align}
excludes the already known solution from further search.
If problem
(\ref{eq:log_min_d3})-(\ref{eq:csat1})-(\ref{eq:csat2})
has no feasible solution, then all optimal solutions of
(\ref{eq:log_min_d3})
have been found. Otherwise, each recently found optimal solution
brings a constraint as (\ref{eq:csat2}), and resolve
(\ref{eq:log_min_d3})-(\ref{eq:csat1})-(\ref{eq:csat2}).
The algorithm stops in a f{\kern0pt}inite number of steps, resulting in
all optimal solutions through binary variables (\ref{eq:log_min_d3}).

Problem (\ref{eq:log_min_alpha2}) can also be rewritten as in (\ref{eq:log_min_d3}):
    \begin{equation}
        \begin{array}{rllllll}
        \min         &&\alpha\\
        \mathrm{s.t.}&& \phi_n(\exp X)\le \alpha, \\
        &&\sum\limits_{i=1}^{n-1} \sum\limits_{j=i+1}^{n} y_{ij}\le K, \\
        && x_{ij}=-x_{ji},&& 1\le i \le j \le n,\\
        && -\bar M\le x_{ij} \le  \bar M, &&1\le i <j \le n,\\
        && -2\bar My_{ij}\le x_{ij}-\bar a_{ij} \le  2\bar My_{ij}, &&1\le i <j \le n,\\
        && y_{ij} \in \{0,1\},&&1\le i <j \le n.\\
        \end{array}
    \label{eq:log_min_alpha3}
    \end{equation}
If  $\phi_n(\exp X)$ is a convex function of the elements (above the main diagonal) of $X$,
then the relaxations of  (\ref{eq:log_min_d3}) and (\ref{eq:log_min_alpha3})
are convex optimization problems, consequently,  (\ref{eq:log_min_d3}) and
(\ref{eq:log_min_alpha3}) are mixed 0-1 convex problems.

\section{Inconsistency index $CR$ of Saaty}

Saaty (1980) proposed to index the inconsistency of pairwise
comparison matrix $A$ of size $n \times n$ by a positive
linear transformation of its largest eigenvalue $\lambda_{\max}$.
The normalized
right eigenvector associated to $\lambda_{\max}$ also plays
an important role, since it provides the estimation of the
weights in the eigenvector method. However, in this paper
weighting methods are not discussed. Saaty (1977) showed
that $\lambda_{\max} \ge n$ and $\lambda_{\max} = n$
if and only if $A$ is consistent. Let us generate a large
number of random pairwise comparison matrices of size $n \times n$, where
each element above the main diagonal are chosen from the
ratio scale $1/9 , 1/8 , 1/7 , \dots, 1/2 , 1, 2, . . . , 8, 9$
with equal probability. Take the largest eigenvalue of each
matrix and let $\overline{\lambda_{\max}}$ denote their
average value. \\
Let $RI_n = (\overline{\lambda_{\max}}-n)/(n-1)$.
Saaty def{\kern0pt}ined the inconsistency of matrix $A$ as
\[
CR_n(A) = \frac{\frac{\lambda_{\max}(A)-n}{n-1}}{RI_n}
\]
being a positive linear transformation of $\lambda_{\max}(A)$.
Then $CR_n(A) \geq 0$ and $CR_n(A) = 0$ if and only if $A$ is consistent.
The heuristic rule of acceptance is $CR_n \leq 0.1$ for all sizes,
also known as the ten percent rule (Saaty, 1980), supported by Vargas' (1982)
statistical analysis. However, some ref{\kern0pt}inements are also known:
$CR_3 \leq 0.05$ for $3\times 3$ matrices
$CR_4 \leq 0.08$ for $4\times 4$ matrices (Saaty, 1994).
Note that any rule of acceptance is somehow heuristic.

Now we apply the results of Section 2 by setting $\phi_n = CR_n$.
Let \linebreak $X\in \log {\cal P}_n$ and let $\lambda_{\max}(\exp X)$ denote
the largest eigenvalue of $A=\exp X$. Then
\begin{align}
\phi_n(\exp X)=\frac {\lambda_{\max}(\exp X)-n}{RI_n(n-1)}. \label{eq:phinexpx}
\end{align}
Boz\'oki et al.\ (2010) showed that $\lambda_{\max}(\exp X)$ is a convex
function of the elements of $X$, therefore, through (\ref{eq:phinexpx}),
$\phi_n(\exp X)$ is a convex function of the elements of $X$, too.

It is proven that (\ref{eq:phinexpx}) implies that both
(\ref{eq:log_min_d3}) and (\ref{eq:log_min_alpha3}) are mixed 0-1 convex
optimization problems. However, they are still challenging from
numerical computational point of view, since $\phi_n(\exp X)$ cannot
be given in an explicit form as $\lambda_{\max}$ values are themselves
computed by iterative methods (Saaty, 1980).
We will show that $\lambda_{\max}$ is not only a limit of an iterative
process, but an optimal solution of a convex optimization problem as well.
The embedded convex optimization problem can be considered together
the embedding optimization problem.

Harker (1987) described the derivatives of $\lambda_{\max}$ with respect to
a matrix element and recommended to change the element with the largest
decrease in $\lambda_{\max}$. The theorems in this section, based on other tools,
can be considered as some extensions of Harker's idea.
Reducing $CR$, being equivalent to decreasing $\lambda_{\max}$, is in the
focus of Xu and Wei (1999) and Cao et al. (2008).

A special case of Frobenius theorem is applied (Saaty, 1977; Sekitani and Yamaki, 1999):

\bigskip
\noindent \textbf{Theorem\ 1.}  \emph{Let $A$ be an $n \times n$ irreducibile nonnegative
matrix and $\lambda_{\max}(A)$ denote the maximal eigenvalue of A. Then the following equalities hold
\begin{align}
\max_{w>0} \min_{i=1,\dots ,n}  \frac{\sum\limits_{j=1}^{n} a_{ij}w_j}{w_i}
 =\lambda_{\max}(A)=  \min_{w>0} \max_{i=1,\dots ,n}
\frac{\sum\limits_{j=1}^{n} a_{ij}w_j}{w_i}. \label{eq:Frobenius}
\end{align}}

\bigskip Since the pairwise comparison matrices are positive, Theorem 1 can be applied.
In order to rewrite the right-hand side of (\ref{eq:Frobenius}), $\bar a_{ij}=\log a_{ij},$
\linebreak $i,j=1,\dots ,n$, and  $z_i=\log w_i,~i=1,\dots ,n$ are used:
\begin{align}
\lambda_{\max}(A)=  \min_{z} \max_{i=1,\dots ,n} \sum\limits_{j=1}^{n} e^{\bar
a_{ij}+z_j-z_i} \label{eq:Frobenius2}
\end{align}

The sum of convex exponential functions in the right-hand side (\ref{eq:Frobenius2}), furthermore,
their maximum are also convex.
Thus, $\lambda_{\max}$ can be determined as the optimum value of a convex optimization problem,
and the form (\ref{eq:Frobenius2}) is equivalent to the optimization problem
\begin{align}
  \min {~\lambda} ~~~\text{s.t}.~~~\sum_{j=1}^n e^{\bar a_{ij}+ z_j-z_i} \leq \lambda, ~ i=1, \dots,n,
   \label{eq:Frobenius3}
\end{align}
where $\lambda$ and $z_{i}, i=1,~\dots,n$ are variables.

Let $\alpha_n$ be given as a threshold of inconsistency index $\phi_n=CR_n$.
Then the constraint
\begin{align}
\phi_n(\exp X)\le \alpha_n \label{eq:phinexpx2a}
\end{align}
from problem (\ref{eq:log_min_d3})
can be transformed by using (\ref{eq:phinexpx}) as
\begin{align}
  \lambda_{\max}(\exp X)\le n+RI_n(n-1)\alpha_n.
   \label{eq:lambdamax}
\end{align}
Denote $\alpha_n^*=n+RI_n(n-1)\alpha_n$. Hence, the formula (\ref{eq:Frobenius2}),
substituting \linebreak $x_{ij}$ = $\bar a_{ij}$, implies an equivalent form
\begin{align}
 \sum\limits_{j=1}^{n} e^{x_{ij}+z_{j}-z_{i}}\le \alpha_n^*,\ i=1,\dots,n.
   \label{eq:lambdamax2}
\end{align}

Let us replace formula (\ref{eq:phinexpx2a})  by (\ref{eq:lambdamax2}) in problem (\ref{eq:log_min_d3}).
We get a mixed 0-1 convex programming problem:
 \begin{equation}
        \begin{array}{rllllll}
        \min         &&\sum\limits_{i=1}^{n-1} \sum\limits_{j=i+1}^{n} y_{ij}\\
        \mathrm{s.t.}&& \sum\limits_{j=1}^{n} e^{x_{ij}+z_{j}-z_{i}}\le \alpha_n^*,&& i=1,\dots,n, \\
        && x_{ij}=-x_{ji},&& 1\le i \le j \le n,\\
        && -\bar M\le x_{ij} \le  \bar M, &&1\le i <j \le n,\\
        && -2\bar My_{ij}\le x_{ij}-\bar a_{ij} \le  2\bar My_{ij}, &&1\le i <j \le n,\\
        && y_{ij} \in \{0,1\},&&1\le i <j \le n.\\
        \end{array}
    \label{eq:log_min_d4}
 \end{equation}

\bigskip
\noindent \textbf{Theorem\ 2.}
\emph{Let $\alpha_n$ denote the acceptance threshold of inconsistency and let \linebreak
$\alpha_n^*=n+RI_n(n-1)\alpha_n$.
Then the optimum value of (\ref{eq:log_min_d4}) gives
the minimal number of the elements to be modif{\kern0pt}ied above the main diagonal in $A$
(and their reciprocals) in order to achieve that $CR_n \leq \alpha_n$.}

\bigskip
Problem (\ref{eq:log_min_alpha3}) can also be rewritten in case of $\phi_n=CR_n$.
In the light of (\ref{eq:phinexpx}), the minimization of $\phi_n$ is equivalent to the minimization
of $\lambda_{\max}$.
Furthermore, program (\ref{eq:Frobenius3}) depending on $\lambda_{\max}$ is used to obtain:

   \begin{equation}
        \begin{array}{rllllll}
        \min         &&\lambda\\
        \mathrm{s.t.}&& \sum\limits_{j=1}^{n} e^{x_{ij}+z_{j}-z_{i}}\le \lambda,&& i=1,\dots,n,  \\
        &&\sum\limits_{i=1}^{n-1} \sum\limits_{j=i+1}^{n} y_{ij}\le K, \\
        && x_{ij}=-x_{ji},&& 1\le i \le j \le n,\\
        && -\bar M\le x_{ij} \le  \bar M, &&1\le i <j \le n,\\
        && -2\bar My_{ij}\le x_{ij}-\bar a_{ij} \le  2\bar My_{ij}, &&1\le i <j \le n,\\
        && y_{ij} \in \{0,1\},&&1\le i <j \le n.\\
        \end{array}
    \label{eq:log_min_alpha4}
    \end{equation}

\bigskip
\noindent \textbf{Theorem \ 3.} \emph{Denote the optimum value of
(\ref{eq:log_min_alpha4}) by $\lambda_{\mathrm{opt}}$, and let \linebreak
$\alpha_{\mathrm{opt}}=\frac{\lambda_{\mathrm{opt}}-n}{RI_n(n-1)}$.
Then $\alpha_{\mathrm{opt}}$ is the minimal value of inconsistency $CR_n$
which can be obtained by the modif{\kern0pt}ication of at most $K$ elements
above the main diagonal of $A$ (and their reciprocals).}

\section{Inconsistency index $CM$ of Koczkodaj and Duszak}

The inconsistency index introduced by Koczkodaj (1993) is based on
$3\times 3 $ submatrices, called \emph{triads}. For the $3\times 3$ pairwise comparison matrix
\begin{align} \left(
\begin{array}{ccc}
    1 & a & b \\
    1/a& 1 & c \\
    1/b & 1/c & 1 \\
\end{array}\label{eq:triad}
             \right)
\end{align}
let
$$
 CM(a,b,c)=\min\left\{\frac{1}{a} \left|a-\frac{b}{c}\right|,\frac{1}{b}\left|b-ac\right|,\frac{1}{c}\left|c-\frac{b}{a}\right|\right\}.
$$
$CM$ can be extended to larger sizes (Duszak and Koczkodaj, 1994):
\begin{align}
    CM(A)=\max\left\lbrace  CM(a_{ij}, a_{ik}, a_{jk}) |~1\le  i< j < k\le n \right\rbrace .
    \label{eq:cma}
\end{align}
Unlike $CR_n$, the construction above does not contain any parameter depending on $n$,
so we dispense with the use of the notation $CM_n$. It is easy to see that $CM$ is an inconsistency index since
$CM(A)\ge 0$ for any $A\in {\cal P}_n$, and $CM(A)=0$ if and only if $A$ is consistent.

For a general triad $(a,b,c)$ let
\begin{align} T(a,b,c)=
\max\left\lbrace \frac{ac}{b}, \frac{b}{ac} \right\rbrace .\label{eq:tabc}
\end{align}
It can be shown (Boz\'oki and Rapcs\'ak, 2008) that there exists a direct relation between
$CM$ and $T$:
\begin{align}
CM(a,b,c)=1-\frac{1}{T(a,b,c)}, \ \
T(a,b,c)=\frac{1}{1-CM(a,b,c)}\label{eq:cmt}.
\end{align}
Since $T(a,b,c)\ge 1$, we get $0\le CM(a,b,c)< 1$, so  $0\le CM(A)< 1$.

Let $(\bar a,\bar b,\bar c)$ denote the logarithmized values of the triad $(a,b,c)$, and let
$$\bar T(\bar a, \bar b,\bar c)=  \max\left\lbrace \bar a+\bar c-\bar b,~
-(\bar a+\bar c-\bar b) \right\rbrace .$$ Then
\begin{align}
T(a,b,c)&= \exp (\bar T(\bar a, \bar b,\bar c) ),\label{eq:tabc2}\\
CM(a,b,c)&=1-\frac{1}{\exp  (\bar T(\bar a, \bar b,\bar c))  }\label{eq:cmabc}.
\end{align}

It is easy to check that even for triads, $CM$ is not a convex function of the logarithmized matrix elements, thus, if we choose the inconsistency index $\phi_n=CM$, then $\phi_n(\exp X)$ appearing in (\ref{eq:log_min_d3}) and (\ref{eq:log_min_alpha3})
is not a convex function of the element of matrix $X$. We show however that by using the univariate function
\begin{align}
f(t)=\frac {1}{1-t}\label{eq:f}
\end{align}
being strictly monotone increasing on the interval $(-\infty,1)$, $f(\phi_n(\exp X))=f(CM(\exp X))$ is already a convex function of the elements of matrix $X$. Then we can change the constraint $$\phi_n(\exp X)\le \alpha_n$$ of problem  (\ref{eq:log_min_d3}) to the convex constraint $$f(\phi_n(\exp X))\le f(\alpha_n).$$
Also, instead of function $\phi_n(\exp
X)$ appearing in problem (\ref{eq:log_min_alpha3}) we can write  $f(\phi_n(\exp X))$ directly, and the value $f^{-1}(\alpha^*)$
computed from the optimal value $\alpha^*$ of the modif{\kern0pt}ied problem is the optimal value of the original problem (\ref{eq:log_min_alpha3}).

To show the  statement above, extend the index $T$ def{\kern0pt}ined in (\ref{eq:tabc}) for arbitrary $n \times n$ pairwise comparison matrix $A$:
\begin{align}
    T(A)=\max\left\lbrace  T(a_{ij}, a_{ik}, a_{jk}) |~1\le  i< j < k\le n \right\rbrace .
    \label{eq:ta}
\end{align}
According to (\ref{eq:cmt}), used there for triads, there is a strictly monotone increasing functional relationship  between $CM$ and $T$. Consequently,
\begin{align}
CM(A)=1-\frac{1}{T(A)}=f^{-1}(T(A)), \ \
T(A)=\frac{1}{1-CM(A)}=f(CM(A))\label{eq:cmta},
\end{align}
where $f$ is the function def{\kern0pt}ined in (\ref{eq:f}).

By expressing $T$ in the logarithmized space, we get
\begin{align}
    T(\exp X)=\max\left\lbrace \max\{ e^{x_{ij}+x_{jk}+x_{ki}}, e^{-x_{ij}-x_{jk}-x_{ki}}\}
     \mid 1\le  i< j < k\le n \right\rbrace .
    \label{eq:texpx}
\end{align}
Since on the right-hand-side of (\ref{eq:texpx}) the maximum of convex functions is taken, $T(\exp X)$ is  convex function of the elements of matrix $X$. Consequently, if we choose the inconsistency index
$\phi_n=CM$, then $f(\phi_n(\exp X))$ is already a convex function, and the problems (\ref{eq:log_min_d3})
and (\ref{eq:log_min_alpha3}) modif{\kern0pt}ied as shown above are already convex mixed 0-1 optimization problems.

Although $CM(\exp X)$ is not convex, it is quasiconvex. To prove it, we show that the lower level sets of $CM(\exp X)$ are convex.
Let $\beta\in [0,1)$ an arbitrary possible value of  $CM(\exp X)$.
Since $f$ is strictly monotone increasing, we have
$$\{X\in \mathbb{R}^{n\times n}\mid CM(\exp X)\le \beta\}=\{X\in \mathbb{R}^{n\times n}\mid f(CM(\exp X))\le f(\beta )\}.$$
Due to the convexity of $T(\exp X)=f(CM(\exp X))$ the above level set are convex, and this implies the quasiconvexity of
$CM(\exp X)$.

\bigskip
\noindent \textbf{Theorem\ 4.}  \emph{$CM(\exp X)$ is quasiconvex on the set of the $n \times
n$ matrices, and \linebreak $T(\exp X)=f(CM(\exp X))$ is convex, where $f$ is def{\kern0pt}ined in (\ref{eq:f}).}

\bigskip

In the following we show that problems (\ref{eq:log_min_d3}) and
(\ref{eq:log_min_alpha3}) can be solved in an easier way, namely,
by solving appropriate linear mixed 0-1 optimization problems.
By exploiting the strictly monotone increasing property of the exponential function,
(\ref{eq:texpx}) can also be written in the following form:
\begin{align}
    T(\exp X)= e^{\max\left\lbrace \max\{{x_{ij}+x_{jk}+x_{ki}}, {-x_{ij}-x_{jk}-x_{ik}}\}
     \mid 1\le  i< j < k\le n \right\rbrace }.
    \label{eq:texpx2}
\end{align}
Now, (\ref{eq:texpx2}) also means that $CM(A)$ can be obtained by determining
the maximum of linear expressions of the elements of matrix  $\bar A=\log A$ and
by applying the exponential function and function $f$ once.

\bigskip
\noindent \textbf{Theorem\ 5.} (Boz\'oki et al.\ 2011a) \emph{ For any $n \times n$ pairwise
comparison matrix $A$, inconsistency index $CM$  can be obtained from the
optimal solution of the following univariate linear program:
    \begin{equation}
        \begin{array}{rllllll}
            \min         && z  \\
            \mathrm{s.t.}&& \bar a_{ij}+\bar a_{jk}+\bar a_{ki}\le z,&&1\le i <j<k \le n, \ \\
                && -(\bar a_{ij}+\bar a_{jk}+\bar a_{ki})\le z &&1\le i <j<k \le n. \\
        \end{array} \label{eq:zopt}
    \end{equation}
Let $z_{\mathrm{opt}}$ be the optimal value of (\ref{eq:zopt}).
Then $CM(A)=1-\frac{1}{\exp(z_{\mathrm{opt}})}$.}

\bigskip

In the following let $\alpha_n$ denote the acceptance threshold associated with the inconsistency index $\phi_n=CM$, and let
\begin{align}
 \alpha_n^*=\log\left(  \frac{1}{1-\alpha_n}\right).
\end{align}

Consider the linear mixed 0-1 optimization problem
\begin{equation}
        \begin{array}{rllllll}
        \min         &&\sum\limits_{i=1}^{n-1} \sum\limits_{j=i+1}^{n} y_{ij}\\
        \mathrm{s.t.}&& x_{ij}+x_{jk}+x_{ki}\le \alpha_n^*,&&1\le i <j<k \le n,\\
            && -(x_{ij}+x_{jk}+x_{ki})\le  \alpha_n^* ,&&1\le i <j<k \le n,\\
        && x_{ij}=-x_{ji},&&1\le i \le j \le n,\\
        && -\bar M\le x_{ij} \le \bar M,&&1\le i <j \le n,\\
        && -2\bar My_{ij}\le x_{ij}-\bar a_{ij} \le 2\bar My_{ij},&&1\le i <j \le n,\\
        && y_{ij} \in \{0,1\},&&1\le i <j \le n.\\
        \end{array}
    \label{eq:mip2}
 \end{equation}

Based on the f{\kern0pt}indings above, the following two theorems follow.

\bigskip
\noindent \textbf{Theorem\ 6.}
\emph{Let $\alpha_n$ denote the acceptance threshold of inconsistency and let \linebreak
$\alpha_n^*=\log(\frac{1}{1-\alpha_n})$.
Then the optimum value of (\ref{eq:mip2}) gives
the minimal number of the elements to be modif{\kern0pt}ied above the main diagonal in $A$
(and their reciprocals) in order to achieve that $CM \leq \alpha_n$.}\\

By some alterations in (\ref{eq:mip2}), the following linear mixed 0-1 optimization problem can be written:
\begin{equation}
        \begin{array}{rllllll}
        \min         && \alpha\\
        \mathrm{s.t.}&& x_{ij}+x_{jk}+x_{ki}\le \alpha,&&1\le i <j<k \le n,\\
            && -(x_{ij}+x_{jk}+x_{ki})\le  \alpha ,&&1\le i <j<k \le n,\\
        && \sum\limits_{i=1}^{n-1} \sum\limits_{j=i+1}^{n} y_{ij}\le K,\\
        && x_{ij}=-x_{ji},&&1\le i \le j \le n,\\
        && -\bar M\le x_{ij} \le \bar M,&&1\le i <j \le n,\\
        && -2\bar My_{ij}\le x_{ij}-\bar a_{ij} \le 2\bar My_{ij},&&1\le i <j \le n,\\
        && y_{ij} \in \{0,1\},&&1\le i <j \le n.\\
        \end{array}
    \label{eq:mip3}
 \end{equation}

\bigskip
\noindent \textbf{Theorem 7.} \emph{Let $\alpha_{\mathrm{opt}}$ denote the
optimum value of (\ref{eq:mip3}).
Then  $1-\frac{1}{\exp(\alpha_{\mathrm{opt}})}$ is
the minimal value of inconsistency $CM$  which
can be obtained by the modif{\kern0pt}ication of at most $K$ elements
above the main diagonal of $A$ (and their reciprocals).}

\section{Inconsistency index $CI$ of Pel\'aez and Lamata}

Similarly to $CM$, the inconsistency index $CI$ proposed by Pel\'aez and Lamata (2003) is also based on triads of form (\ref{eq:triad}). It is easy to see that the determinant of the triad (\ref{eq:triad}) is nonnegative, and it is zero if and only if the triad is consistent. Based on this interesting property, Pel\'aez and Lamata (2003) proposed to characterize the inconsistency of a pairwise comparison matrix $A \in {\cal P}_n$ by the average of the determinants of the triads of matrix $A$:
\begin{equation}
CI_n(A) =
\begin{cases}
\det(A), &  \textrm{for} ~n=3, \\
\frac{1}{NT(n)} \sum\limits_{i=1}^{NT(n)} \det(\Gamma_i), &  \textrm{for} ~n>3,\\
\end{cases}
\label{eq:CI_def}
\end{equation}
where $\Gamma_i$, $i=1,\dots,NT(n)$ denote the triads of matrix $A$, and $NT(n)={n \choose 3}$ is the number of triads in $A$.

We show that $CI$ is a convex function of the logarithmized matrix elements, thus if  the inconsistency index $\phi_n=CI_n$ is chosen, then $\phi_n(\exp X)$ appearing in problems  (\ref{eq:log_min_d3}) and (\ref{eq:log_min_alpha3}) is a convex function of the elements of matrix $X$.

The determinant of triad  $\Gamma \in {\cal P}_3$ comparing objects $(i,j,k)$ can be written as
\begin{equation}
 \det(\Gamma)=\frac{a_{ik}}{a_{ij}a_{jk}} + \frac{a_{ij} a_{jk}}{a_{ik}}-2. \label{eq:sarrus}
\end{equation}

Let $ X= \log \Gamma \in\log ~{\cal P}_3$, i.e.,  $\Gamma= \exp X$. Equation (\ref{eq:sarrus}) can be reformulated as a convex function of the elements of $X$:
\begin{equation}
 \det(\exp X)=e^{x_{ik} - x_{ij} -  x_{jk}} + e^{x_{ij} + x_{jk} - x_{ik}}-2. \label{eq:sarrus_log}
\end{equation}

Let $\alpha_n$ be a given acceptance threshold for the inconsistency index \linebreak $\phi_n=CI_n$.
According to (\ref{eq:CI_def}) and (\ref{eq:sarrus_log}), the constraint
\begin{equation}
\phi_n(\exp X)\le \alpha_n \label{eq:phinexpx2}
\end{equation}
appearing in (\ref{eq:log_min_d3}) can be expressed as
\begin{equation}
  \frac{1}{{n \choose 3}} \sum\limits_{i=1}^{n-2}\sum\limits_{j=i+1}^{n-1}  \sum\limits_{k=j+1}^{n}  \left( e^{x_{ik}-x_{ij}-x_{jk}} + e^{x_{ij}+x_{jk}-x_{ik}}-2 \right) \le \alpha_n.
   \label{eq:CI_eq_log}
\end{equation}
By denoting $\alpha_n^*=(\alpha_n+2) {n \choose 3}$, (\ref{eq:CI_eq_log}) can be simplif{\kern0pt}ied as
\begin{equation}
  \sum\limits_{i=1}^{n-2}\sum\limits_{j=i+1}^{n-1}  \sum\limits_{k=j+1}^{n}  \left( e^{x_{ik}-x_{ij}-x_{jk}} + e^{x_{ij}+x_{jk}-x_{ik}} \right) \le \alpha_n^*,
   \label{eq:CI_eq_log_2}
\end{equation}
and inserting it into (\ref{eq:log_min_d3}), we get the  mixed 0-1 convex optimization problem
%\begin{small}
    \begin{equation}
        \begin{array}{lllll}
        \min        &\sum\limits_{i=1}^{n-1} \sum\limits_{j=i+1}^{n} y_{ij}\\
        \mathrm{s.t.}& \sum\limits_{i=1}^{n-2}\sum\limits_{j=i+1}^{n-1}  \sum\limits_{k=j+1}^{n}  \left( e^{x_{ik}-x_{ij}-x_{jk}} + e^{x_{ij}+x_{jk}-x_{ik}} \right) \le \alpha_n^* ,\\
        & x_{ij}=-x_{ji},\ &\hskip -3cm 1\le i <j \le n,\\
        & -\bar M\le x_{ij} \le  \bar M,\ &\hskip -3cm 1\le i <j \le n,\\
        & -2\bar My_{ij}\le x_{ij}-\bar a_{ij} \le  2\bar My_{ij}, &\hskip -3cm 1\le i <j \le n,\\
        & y_{ij} \in \{0,1\},&\hskip -3cm 1\le i <j \le n.\\
        \end{array}
    \label{eq:minlp1_CI}
\end{equation}
%\end{small}

\bigskip
\noindent \textbf{Theorem\ 8.}
\emph{Let $\alpha_n$ denote the acceptance threshold of inconsistency and let \linebreak
$\alpha_n^*=(\alpha_n+2) {n \choose 3}$.
Then the optimum value of (\ref{eq:minlp1_CI}) gives
the minimal number of the elements to be modif{\kern0pt}ied above the main diagonal in $A$
(and their reciprocals) in order to achieve that
$CI \leq \alpha_n$. }

\bigskip
In the same way as for other inconsistency indices, the following  mixed 0-1 convex optimization problem can also be considered:
 \begin{equation}
        \begin{array}{lllll}
        \min         && \alpha\\
        \mathrm{s.t.}&& \sum\limits_{i=1}^{n-2}\sum\limits_{j=i+1}^{n-1}  \sum\limits_{k=j+1}^{n}  \left( e^{x_{ik}-x_{ij}-x_{jk}} + e^{x_{ij}+x_{jk}-x_{ik}} \right) \le \alpha ,\\
                    &&  x_{ij}=-x_{ji}, &  \!\!\!\!\!\!\!\!\!\!\!\!\!\!\!\!\!\!\!\!\!\!\!\!\!\!\!\!\!\!\!\!\!\!\!\!\!\! 1\le i <j \le n,\\
        && -\bar M\le x_{ij} \le  \bar M,\ & \!\!\!\!\!\!\!\!\!\!\!\!\!\!\!\!\!\!\!\!\!\!\!\!\!\!\!\!\!\!\!\!\!\!\!\!\!\!  1\le i <j \le n,\\
        && -2\bar My_{ij}\le x_{ij}-\bar a_{ij} \le  2\bar My_{ij}, & \!\!\!\!\!\!\!\!\!\!\!\!\!\!\!\!\!\!\!\!\!\!\!\!\!\!\!\!\!\!\!\!\!\!\!\!\!\! 1\le i <j \le n,\\
        && y_{ij} \in \{0,1\}, & \!\!\!\!\!\!\!\!\!\!\!\!\!\!\!\!\!\!\!\!\!\!\!\!\!\!\!\!\!\!\!\!\!\!\!\!\!\!  1\le i <j \le n,\\
        && \sum\limits_{i=1}^{n-1} \sum\limits_{j=i+1}^{n} y_{ij} \le K.
        \end{array}
    \label{eq:minlp2_CI}
\end{equation}

\bigskip
\noindent \textbf{Theorem 9.} \emph{Let $\alpha_{\mathrm{opt}}$ denote the optimum
value of (\ref{eq:minlp2_CI}). Then  $ \frac{\alpha_{\mathrm{opt}}}{{n \choose 3}}   -2$
is the minimal value of inconsistency $CI$  which
can be obtained by the modif{\kern0pt}ication of at most $K$ elements
above the main diagonal of $A$ (and their reciprocals).}

\section{A numerical example}

Our approach is also presented on a classic numerical example from the book of Saaty (1980),
for the inconsistency index $CR$. Table 1 contains pairwise comparison values of six cities concerning their distances from Philadelphia.
As an example, the evaluator judged that the distance between London and Philadelphia is f{\kern0pt}ive times greater than that between Chicago and Philadelphia.

\bigskip
\bigskip
\noindent \textbf{Table 1.} Comparison of distances of cities from Philadelphia
\newline
\begin{tabular}{lllllll  }
 \hline
 & Cairo&Tokyo&Chicago&San Francisco&London&Montreal\\
 \hline
 Cairo  &1 &1/3 & 8 & 3&3&7 \\
 Tokyo& 3& 1& 9 &3 &3&9 \\
 Chicago & 1/8&1/9 & 1 &1/6 &1/5&2\\
 San Francisco  & 1/3&1/3 &6  &1 &1/3&6\\
 London  & 1/3& 1/3& 5 & 3&1& 6\\
 Montreal& 1/7& 1/9& 1/2&1/6  &1/6 &1 \\
 \hline
 \end{tabular}

 \bigskip
\bigskip
Let $A$ denote the pairwise comparison matrix concerning Table\ 1. We get that $\lambda_{\max}(A)=6.4536$, and from $RI_6=1.24$, also $CR(A)=0.0732$. Since the value of
$CR(A)$  is signif{\kern0pt}icantly below the $10\%$ threshold, we can consider the inconsistency of $A$ acceptable.

Let $A^{(1)}$ denote the matrix obtained from $A$ by exchanging the elements
$a_{1,2}$ (and $a_{2,1}$). This is a typical mistake at f{\kern0pt}illing-in a pairwise comparison matrix. For the matrix $A^{(1)}$, we get $CR(A^{(1)})=0.0811$. Therefore, although the level of inconsistency of $A^{(1)}$ has increased as consequence of the data-recording error, it is still below the acceptance level of 10\%. In this case the proposed methodology is not able to detect the mistake, and $A^{(1)}$ is still accepted.

Consider now the case when $a_{1,3}$ and $a_{3,1}$ are exchanged, say  by accident, in the matrix $A$. Let $A^{(2)}$ denote the matrix obtained in this way. Then $CR(A^{(2)})=0.5800$, which is well over the acceptance level of 10\%, and it refers to a rough inconsistency in the matrix. By solving the corresponding problem (\ref{eq:log_min_d4}), we obtain that
the inconsistency of $A^{(2)}$ can be pushed below the critical 10\% by modifying a single element (and its reciprocal).  This element is just in the spoilt position $a_{1,3}$.  It can also be shown that this is the single optimal solution to problem (\ref{eq:log_min_d4}) considering the 0-1 variables. Consequently, the proposed methodology has detected the single possible element for the case of correcting in a single position (and in its reciprocal). It also turned out that this single position is just the one  of the values exchanged by accident.

In the previous example the spoilt matrix caused a rough increase of the inconsistency. In this view, it is not surprising that the proposed method of{\kern0pt}fers a unique way of repairing. However, at smaller increase of inconsistency the situation is not that obvious.

Assume now that the element $a_{1,3}$ of matrix $A$ is changed to 2 instead of the value 1/8 of the previous example. This is a smaller dif{\kern0pt}ference in relation to the original value 8, the increase of the inconsistency of the modif{\kern0pt}ied matrix, denoted by $A^{(3)}$, is also less: $CR(A^{(3)})=0.1078$.
The inconsistency of $A^{(3)}$ barely exceeds the critical level 10\%, therefore, one would expect that by the modif{\kern0pt}ication of a single element can make the inconsistency decrease below 10\%, and also that several positions are eligible for this purpose.
Indeed, the optimal value of the relating problem (\ref{eq:log_min_d4}) is 1,
and by resolving the problem after adding the constraints (\ref{eq:csat1}) and (\ref{eq:csat2}) we f{\kern0pt}ind that problem  (\ref{eq:log_min_d4}) has 6 dif{\kern0pt}ferent optimal solutions according to the binary variables.
Namely, the inconsistency of matrix $A^{(3)}$ decreases below 10\% not only by modifying $a_{1,3}$, but also by modifying any single element of $\{a_{1,4}, a_{1,5}, a_{2,6}, a_{3,4}, a_{4,5}\}$.
In the ideal case, the evaluator spots the data-recording error in position $a_{1,3}$ immediately. If not, then s/he may have to reconsider the evaluation of each of the 6 positions, but it is still fewer than the 15 possible positions in the upper triangular part of the matrix.

\bigskip
\bigskip

\begin{hangparas}{.3in}{1}

\textbf{References}

\bigskip

%\bibitem{BozokiFulopKoczkodaj2011}
Boz\'oki S, F\"ul\"op J, Koczkodaj WW (2011a) LP-based consistency-driven
supervision for incomplete pairwise comparison matrices.
Mathematical and Computer Modelling 54(1-2):789--793

Boz\'oki S, F\"ul\"op J, Poesz A (2011b)
On pairwise comparison matrices that can be made consistent by the modif{\kern0pt}ication
of a few elements. Central European Journal of Operations Research 19(2):157--175

Boz\'oki S, F\"ul\"op J, Poesz A (2012)
Convexity properties related to pairwise comparison matrices
of acceptable inconsistency and applications, \linebreak
(in Hungarian,
Elfogadhat\'o inkonzisztenci\'aj\'u p\'aros \"osszehasonl\'it\'as m\'at-rixokkal
kapcsolatos konvexit\'asi tulajdons\'agok \'es azok alkalmaz\'asai).
In: Solymosi T, Temesi J (eds.) Egyens\'uly \'es optimum:
Tanulm\'anyok Forg\'o Ferenc 70.~sz\"ulet\'esnapj\'ara, Aula Kiad\'o, pp.~169--184

%\bibitem{BoFuRo2010}
Boz\'oki S, F\"ul\"op J, R\'onyai L (2010)
On optimal completions of incomplete pairwise comparison matrices.
Mathematical and Computer Modelling 52(1-2):318--333

%\bibitem{BoRa08}
Boz\'oki S, Rapcs\'ak T (2008) On Saaty's and Koczkodaj's
inconsistencies of pairwise comparison matrices.
Journal of Global Optimization 42(2):157--175

%\bibitem{BrunelliFedrizzi2011}
Brunelli M, Fedrizzi M (2011) Characterizing properties for inconsistency indices
in the AHP. Proceedings of the 11th International Symposium on the AHP,
Sorrento, Naples, Italy, June 15-18, 2011
%http://204.202.238.22/isahp2011/dati/pdf/55\_064\_Fedrizzi.pdf

Brunelli M, Fedrizzi M (2013a)
Axiomatic properties of inconsistency indices for pairwise comparisons.
Submitted, http://arxiv.org/abs/1306.6852

Brunelli M, Canal L, Fedrizzi M (2013b)
Inconsistency indices for pairwise comparison matrices:
a numerical study.
Annals of Operations Research, published online f{\kern0pt}irst,
DOI 10.1007/s10479-013-1329-0

Cao D, Leung LC, Law JS (2008)
Modifying inconsistent comparison matrix in analytic hierarchy process: A heuristic approach.
Decision Support Systems 44(4):944--953

%\bibitem{Chu1997}
Chu MT (1997) On the optimal consistent approximation to pairwise comparison
matrices.
Linear Algebra and its Applications 272(1-3):155--168

%\bibitem{DuKo94}
Duszak Z, Koczkodaj WW (1994) Generalization of a new def{\kern0pt}inition of consistency
for pairwise comparisons.
Information Processing Letters 52(5):273--276

Harker PT (1987)
Derivatives of the Perron root of a positive reciprocal matrix: With application to the analytic hierarchy process.
Applied Mathematics and Computation 22(2-3):217--232

%\bibitem{Koczkodaj93}
Koczkodaj WW (1993) A new def{\kern0pt}inition of consistency of pairwise
comparisons.
Mathematical and Computer Modelling 18(7):79--84

%\bibitem{Pelaez}
Pel\'aez JI, Lamata MT (2003) A new measure of consistency for positive reciprocal matrices.
Computers and Mathematics with Applications 46:1839--1845

Saaty TL (1977)
A scaling method for priorities in hierarchical structures.
Journal of Mathematical Psychology 15(3):234--281

%\bibitem{Saaty1980}
Saaty TL (1980) The Analytic Hierarchy Process, McGraw-Hill, New York

%\bibitem{Saaty1994}
Saaty TL (1994) Fundamentals of Decision Making, RSW Publications, Pittsburgh

%\bibitem{SekitaniYamaki1999}
Sekitani K, Yamaki N (1999)  A logical interpretation for the eigenvalue method
in AHP.
Journal of the Operations Research Society of Japan 42(2):219--232

%\bibitem{Vargas1982}
Vargas LG (1982) Reciprocal matrices with random coef{\kern0pt}f{\kern0pt}icients.
Mathematical Modelling 3(1):69--81

Xu Z, Wei C (1999)
A consistency improving method in the analytic hierarchy process.
European Journal of Operational Research
116(2):443--449

\end{hangparas}

\end{document}